\documentclass[a4paper,11pt,reqno]{amsart}

\newtheorem{theorem}[equation]{Theorem}

\newtheorem{corollary}[equation]{Corollary}
\newtheorem{proposition}[equation]{Proposition}
\numberwithin{equation}{section}

\theoremstyle{definition}
\newtheorem{definition}[equation]{Definition}

\newtheorem{remark}[equation]{Remark}

\newcommand{\cV}{{\mathcal V}}
\newcommand{\g}{{\mathfrak g}}

 \DeclareMathOperator{\Der}{Der}
 \DeclareMathOperator{\End}{End}
 
 \DeclareMathOperator{\Mat}{Mat}

\newenvironment{romanenumerate}
 {\begin{enumerate}
 }{\end{enumerate}}

\newenvironment{romanprimeenumerate}
 {\begin{enumerate}
 }{\end{enumerate}}

\begin{document}

\title[$(-1,-1)$-balanced Freudenthal Kantor triple systems]{$(-1,-1)$-Balanced
Freudenthal Kantor triple systems and noncommutative Jordan
algebras}

\author[A.~Elduque]{Alberto Elduque$^{\star}$}
 \thanks{$^{\star}$ Supported by the Spanish Ministerio de Ciencia y
Tecnolog\'{\i}a and FEDER (BFM 2001-3239-C03-03)}
 \address{Departamento de Matem\'aticas, Universidad de
Zaragoza, 50009 Zaragoza, Spain}
 \email{elduque@unizar.es}

\author[N.~Kamiya]{Noriaki Kamiya$^{\circ}$}
 \thanks{$^{\circ}$ Partially supported by a Grant-in-Aid for Science Research
 (no.~15540037(C)(2)) from the Japan Society for the Promotion of Science}
 \address{Center for Mathematical Science, The University of Aizu.
 Aizu-Wakamatsu, 965-8580 Japan}
  \email{kamiya@u-aizu.ac.jp}

\author[S.~Okubo]{Susumu Okubo$^{\ast}$}
 \thanks{$^{\ast}$ Supported in part by U.S.~Department of Energy Grant No.
 DE-FG02-91ER40685}
 \address{Department of Physics and Astronomy, University of
 Rochester, Rochester, NY 14627, USA}
 \email{okubo@pas.rochester.edu}

\date{\today}



\begin{abstract}
A noncommutative Jordan algebra of a specific type is attached to
any $(-1,-1)$-balanced Freudenthal Kantor triple system, in such a
way that the triple product in this system is determined by the
binary product in the algebra. Over fields of characteristic zero,
the simple noncommutative Jordan algebras of this type are
classified.
\end{abstract}

\maketitle


\section{Introduction}

The well-known Tits-Kantor-Koecher construction \cite{T,K,Ko}
relates Jordan systems to $3$-graded Lie algebras. In
\cite{Kantor}, several models of exceptional Lie algebras with a
$5$-grading
\[
\g=\g_{-2}\oplus \g_{-1}\oplus \g_0\oplus\g_1\oplus \g_2
\]
are given, based on \emph{generalized Jordan triple systems}:

\begin{definition} A vector space $J$ over a field $F$, endowed
with a trilinear operation $J\times J\times J\rightarrow J$,
$(x,y,z)\mapsto xyz$, is said to be a \emph{generalized Jordan
triple system} (GJTS for short) if it satisfies the identity:
\begin{equation}\label{eq:GJTS}
uv(xyz)=(uvx)yz-x(vuy)z+xy(uvz)
\end{equation}
for any $u,v,x,y,z\in J$.
\end{definition}

A (linear) \emph{Jordan triple system} is then a generalized
Jordan triple system  with the added constraint:
\begin{equation}\label{eq:JTS}
xyz=zyx
\end{equation}
for any $x,y,z$.

Unless otherwise stated, all the algebras and algebraic systems
considered will be assumed to be defined over a ground field $F$
of characteristic not $2$.

Given two elements $a,b$ in a GJTS $J$, consider the linear maps
$l_{a,b},k_{a,b}:J\rightarrow J$ given by $l_{a,b}c=abc$,
$k_{a,b}c=acb+bca$. Thus, equation \eqref{eq:GJTS} is equivalent
to
\begin{equation}\label{eq:FKTS1}
[l_{u,v},l_{x,y}]=l_{l_{u,v}x,y}-l_{x,l_{v,u}y}
\end{equation}
for any $u,v,x,y$.

\begin{definition}
A generalized Jordan triple system $J$ is said to be a
\emph{$(-1,-1)$ Freudenthal Kantor triple system} ($(-1,-1)$-FKTS
for short) if it satisfies
\begin{equation}\label{eq:FKTS2}
l_{d,c}k_{a,b}+k_{a,b}l_{c,d}=k_{k_{a,b}c,d}
\end{equation}
for any $a,b,c,d\in J$.
\end{definition}

The more general concept of $(\varepsilon,\delta)$ Freudenthal
Kantor triple system was introduced in \cite{YO} for
$\varepsilon,\delta=\pm 1$.

In case there is a nonzero symmetric bilinear form $\langle .
\vert .\rangle$ defined on a $(-1,-1)$-FKTS  such that
\begin{equation}\label{eq:FKTS3}
k_{a,b}c=\dfrac{1}{2}\langle a\vert b\rangle c
\end{equation}
for any $a,b,c$, the $(-1,-1)$-FKTS is said to be \emph{balanced}.

In this case, equations \eqref{eq:FKTS2} and \eqref{eq:FKTS3}
simplfy to
\begin{equation}\label{eq:FKTSdos}
xxy=xyx=\langle x\vert x\rangle y
\end{equation}
for any $x,y$.

Therefore, a \emph{$(-1,-1)$-balanced Freudenthal Kantor triple
system} (or \linebreak $(-1,-1)$-BFKTS for short) is a vector
spaced over a field $F$ endowed with a trilinear operation $xyz$
and a nonzero symmetric bilinear form $\langle .\vert .\rangle$
such that both \eqref{eq:GJTS} and \eqref{eq:FKTSdos} hold.

Some simple $(-1,-1)$-BFKTS's were used in \cite{KOk} to give
models of the simple complex exceptional Lie superalgebras of
classical type: $D(2,1;\alpha)$, $G(3)$ and $F(4)$. Furthermore,
this close relationship between $(-1,-1)$-BFKTS's and some Lie
superalgebras was used in \cite{EKO} to provide the classification
of the simple finite dimensional $(-1,-1)$-BFKTS's over fields of
characteristic zero.

The aim of this paper is to show that a quadratic noncommutative
Jordan algebra in a specific variety can be attached to any
$(-1,-1)$-BFKTS in such a way that the triple product is
determined by the (binary) multiplication of the algebra. The
classification of the simple finite dimensional quadratic
noncommutative Jordan algebras over fields of characteristic $0$
in this variety will be deduced too from the known classification
of the simple $(-1,-1)$-BFKTS's.

The next section will be devoted to introduce the variety $\cV$ of
noncommutative Jordan algebras that will be relevant for our
purposes. Then in Section 3, the relationship between some GJTS's
and algebras in $\cV$ will be studied, while in Section 4 the
attention will be restricted to $(-1,-1)$-BFKTS's. The last
Section will deal with the classification of the simple finite
dimensional quadratic noncommutative Jordan algebras in the
variety $\cV$ over fields of characteristic $0$.

\section{A variety of noncommutative Jordan algebras}

Given any algebra $A$ over a field $F$ (always of characteristic
$\ne 2$), let $L_x$ and $R_x$ denote the left and right
multiplications by $x$: $L_x(y)=xy$, $R_x(y)=yx$; and let
$(x,y,z)$ denote the \emph{associator} of the elements $x,y,z$:
$(x,y,z)=(xy)z-x(yz)$.

Recall that the algebra $A$ is a \emph{noncommutative Jordan
algebra} if it is flexible:
\begin{equation}\label{eq:flexible}
(x,y,x)=0
\end{equation}
for any $x,y\in A$, and satisfies the Jordan identity:
\begin{equation}\label{eq:Jordan}
(x,y,x^2)=0
\end{equation}
for any $x,y\in A$.

The Jordan identity is equivalent to the condition
$[L_x,R_{x^2}]=0$ for any $x\in A$, while the flexibility amounts
to $[L_x,R_x]=0$ for any $x\in A$. Also, by flexibility
\[
(x,x,y)+(y,x,x)=0
\]
for any $x,y\in A$, or $L_{x^2}-L_x^2=R_{x^2}-R_x^2$. Thus, in the
presence of flexibility, $[L_x,R_{x^2}]=0$ if and only if
\begin{equation}\label{eq:LxLx2}
[L_x,L_{x^2}]=0
\end{equation}
for any $x\in A$.

Given two elements $x,y$ of a flexible algebra $A$, consider the
linear map $A\rightarrow A$ given by:
\begin{equation}\label{eq:Dxy}
\begin{split}
D_{x,y}&=L_{[x,y]}-[L_x,L_y]\\
&=(x,y,\,.\,)-(y,x,\,.\,)\\
&=(\,.\,,x,y)-(\,.\,,y,x)\\
&=-R_{[x,y]}-[R_x,R_y]
\end{split}
\end{equation}
Notice that \eqref{eq:LxLx2} is equivalent to the condition
$D_{x^2,x}=0$.

Let $\cV$ be the variety of those noncommutative Jordan algebras
$A$ over a field $F$ satisfying that
\begin{equation}\label{eq:der}
\text{$D_{x,y}$ is a derivation of $A$ for any $x,y\in A$.}
\end{equation}
This is the variety that will be relevant in what follows.

\begin{theorem}
Let $A$ be an algebra in $\cV$, then for any $x,y,z\in A$,
\begin{equation}\label{eq:cyclic}
D_{xy,z}+D_{yz,x}+D_{zx,y}=0.
\end{equation}
\end{theorem}
\begin{proof}
The fact that $D_{x,y}=L_{[x,y]}-[L_x,L_y]\in \Der A$ (the Lie
algebra of derivations of $A$) for any $x,y\in A$ is equivalent to
the validity of
\[
\left[L_{[x,y]},L_z\right]-\left[[L_x,L_y],L_z\right]
=L_{(x,y,z)-(y,x,z)}
\]
for any $x,y,z\in A$. Permute cyclically $x,y,z$ and add the
resulting equations to get
\begin{equation}\label{eq:Lxyz}
\begin{split}
\left[L_{[x,y]},L_z\right]+\left[L_{[y,z]},L_x\right]&+\left[L_{[z,x]},L_y\right]\\
&= L_{(x,y,z)+(y,z,x)+(z,x,y)-(y,x,z)-(z,y,x)-(x,z,y)} .
\end{split}
\end{equation}
But in any algebra,
\begin{multline*}
(x,y,z)+(y,z,x)+(z,x,y)-(y,x,z)-(z,y,x)-(x,z,y)\\
=[[x,y],z]+[[y,z],x]+[[z,x],y],
\end{multline*}
so \eqref{eq:Lxyz} is equivalent to
\[
\left[L_{[x,y]},L_z\right]+\left[L_{[y,z]},L_x\right]+\left[L_{[z,x]},L_y\right]
= L_{[[x,y],z]+[[y,z],x]+[[z,x],y]} ,
\]
or to
\begin{equation}\label{eq:Dxyz}
D_{[x,y],z}+D_{[y,z],x}+D_{[z,x],y}=0.
\end{equation}
But $D_{x^2,x}=0$ for any $x\in A$, since $A$ is a noncommutative
Jordan algebra, so that, by linearization, $D_{x\circ
y},x+D_{x^2,y}=0$ for any $x,y\in A$, where $x\circ y=xy+yx$, and
\begin{equation}\label{eq:Dxoyz}
D_{x\circ y,z}+D_{y\circ z,x}+D_{z\circ x,y}=0
\end{equation}
for any $x,y,z\in A$.

The result now follows by adding up \eqref{eq:Dxyz} and
\eqref{eq:Dxoyz}.
\end{proof}

The algebras $A$ endowed with a skew symmetric bilinear map $D:
A\times A\rightarrow \Der A$, $(x,y)\mapsto D_{x,y}$ satisfying
\eqref{eq:cyclic} have been named \emph{generalized structurable
algebras} in \cite{Kamiya}. The variety of generalized
structurable algebras includes the most usual varieties of
nonassociative algebras.

\begin{corollary}
Any algebra in $\cV$ is a generalized structurable algebra.
\end{corollary}

\section{Generalized Jordan triple systems}

It will be shown in this section the close connection of the
algebras in the variety $\cV$ with some generalized Jordan triple
systems:

\begin{theorem}\label{GJTS}
Let $J$ be a generalized Jordan triple system over a field $F$ of
characteristic $\ne 2,3$ which contains an element $e\in J$ such
that:
\begin{romanenumerate}
\item $eee=e$,
\item $eex=xee$ for any $x\in J$,
\item the map $U_e: x\mapsto exe$ is onto.
\end{romanenumerate}
Then the homotope algebra $J^{(e)}$, defined on the vector space
$J$ with multiplication given by $x\cdot y=xey$ for any $x,y\in
J$, belongs to the variety $\cV$ and is unital with $1=e$.
Moreover, the map $x\mapsto \bar x=exe$ is an involution of
$J^{(e)}$ and the triple product in $J$ satisfies
\begin{equation}\label{eq:bintriple}
xyz=x\cdot(\bar y\cdot z)-\bar y\cdot (x\cdot z)+(\bar y\cdot
x)\cdot z,
\end{equation}
for any $x,y,z\in J$.

Conversely, let $(A,\cdot)$ be a unital algebra in $\cV$ over a
field $F$ of characteristic $\ne 2$, with unity element $1_A$ and
endowed with an involution $x\mapsto \bar x$, and define a triple
product on $A$ by means of \eqref{eq:bintriple}. Then $A$ becomes
a GJTS and satisfies conditions \textup{(i)--(iii)} above with
$e=1_A$ and $U_ex=\bar x$ for any $x\in A$.
\end{theorem}

\begin{proof}
Let $J$ be a GJTS satisfying the conditions above. Then
\eqref{eq:GJTS}, together with (i) and (ii), give:
\begin{align}
 \begin{aligned}
  ee(exe)&=(eee)xe-e(eex)e+ex(eee)\\
   &=2exe-e(eex)e,
 \end{aligned}\label{eq:G1}\\
 \begin{aligned}
  exe=ex(eee)&=(exe)ee-e(xee)e+ee(exe)\\
    &=2ee(exe)-e(eex)e.
 \end{aligned}\label{eq:G2}
\end{align}
Since the characteristic of $F$ is $\ne 3$, \eqref{eq:G1} and
\eqref{eq:G2} imply
\begin{equation}\label{eq:G3}
exe=ee(exe)=e(eex)e
\end{equation}
for any $x\in J$. Since $U_e:x\mapsto exe$ is onto, we get
\begin{equation}\label{eq:G4}
x=eex\, (=xee),
\end{equation}
so that
\[
x=e\cdot x=x\cdot e
\]
for any $x\in J$ and $e$ is the unity element of the homotope
algebra $J^{(e)}$.

Now, by \eqref{eq:G4} and \eqref{eq:GJTS},
\[
\begin{split}
x=xe(eee)&=(xee)ee-e(exe)e+ee(xee)\\
   &= 2x-e(exe)e
\end{split}
\]
so $\bar{\bar x}=x$ for any $x\in J$, where $\bar x=U_ex=exe$.
Also, for any $x,y\in J$,
\[
\left\{
\begin{aligned}
x\cdot y=xe(eey)&=(xee)ey-e(exe)y+ee(xey)\\
    &= x\cdot y-e\bar x y+x\cdot y,\\
x\cdot y=xe(yee)&= (xey)ee-y(exe)e+ye(xee)\\
   &=x\cdot y-y\bar x e+y\cdot x.
\end{aligned}\right.
\]
Hence, for any $x,y\in J$,
\begin{equation}\label{eq:G5}
\left\{
\begin{aligned}
exy&=\bar x\cdot y,\\
xye&=x\cdot \bar y,
\end{aligned}\right.
\end{equation}
while
\[
\begin{split}
\bar x\cdot \bar y=ex(eye)&=(exe)ye-e(xey)e+ey(exe)\\
&=\bar x ye-\overline{x\cdot y}+ey\bar x\\
&=\bar x\cdot \bar y-\overline{x\cdot y}+\bar y\bar x,
\end{split}
\]
and this shows that $\overline{x\cdot y}=\bar y\cdot \bar x$ for
any $x,y\in J$. Therefore, the map $x\mapsto \bar x$ is an
involution of $J^{(e)}$.

Besides, for any $x,y,z\in J$,
\[
\begin{split}
y\cdot (x\cdot z)=ye(xez)&=(yex)ez-x(eye)z+xe(yez)\\
&=(y\cdot x)\cdot z-x\bar y z+x\cdot (y\cdot z),
\end{split}
\]
so, substituting $y$ by $\bar y$ we obtain \eqref{eq:bintriple}.

But also, because of \eqref{eq:G5},
\[
\begin{split}
z\cdot(x\cdot \bar y)=e\bar z(xye)&=(e\bar zx)ye-
     x(\bar z ey)e+xy(e\bar z e)\\
     &=(z\cdot x)\cdot \bar y -x\cdot \overline{(\bar z\cdot y)}
       +xyz,
\end{split}
\]
so
\begin{equation}\label{eq:G6}
xyz=x\cdot(\bar y\cdot z)-(z\cdot x)\cdot\bar y +z\cdot(x\cdot
\bar y)
\end{equation}
for any $x,y,z$, since $x\mapsto \bar x$ is an involution.
Equations \eqref{eq:bintriple} and \eqref{eq:G6} yield, for any
$x,y,z\in J$,
\begin{equation}\label{eq:G7}
\left(\bar y,x,z\right)^\cdot =-\left(z,x,\bar y\right)^\cdot,
\end{equation}
where $(a,b,c)^\cdot =(a\cdot b)\cdot c-a\cdot (b\cdot c)$ is the
associator in $J^{(e)}$. Identity \eqref{eq:G7} is equivalent to
the flexible law in $J^{(e)}$.

Moreover, \eqref{eq:bintriple} is equivalent to
\begin{equation}\label{eq:G8}
l_{x,y}=[L_x,L_{\bar y}]+L_{\bar y\cdot x}=
 -D_{x,\bar y}+L_{x\cdot\bar y},
\end{equation}
where $L_a$ denotes the left multiplication by $a$ in $J^{(e)}$
and $D_{x,y}=L_{[x,y]}-[L_x,L_y]$ as in \eqref{eq:Dxy}. Notice
that $L_a=l_{a,e}$ and $L_{\bar z}=l_{e,z}$ (because of
\eqref{eq:G5}). Thus,
\begin{equation}\label{eq:G9}
\begin{split}
[D_{x,y},l_{z,t}]
  &=-[l_{x,\bar y},l_{z,t}] + [l_{x\cdot y,e},l_{z,t}]\\
  &=-l_{x\bar y z,t}+l_{z,\bar y xt}+
    l_{(x\cdot y)ez,t}-l_{z,e(x\cdot y)t}\\
  &=l_{-x\bar y z+(x\cdot y)\cdot z,t}
   +l_{z,\bar yxt-(\bar y\cdot\bar x)\cdot t}\\
  &=l_{D_{x,y}z,t}-l_{z,D_{\bar y,\bar x}t}
\end{split}
\end{equation}
for any $x,y,z,t\in J$. With $t=e$, this shows that
$[D_{x,y},L_z]=L_{D_{x,y}z}$, since $D_{a,b}e=0$ for any $a,b\in
J$. Hence $D_{x,y}\in \Der J^{(e)}$ for any $x,y\in J$.

Finally, \eqref{eq:G8} shows that $l_{x,\bar x}=L_{x\cdot x}$, so
\[
\begin{split}
[L_x,L_{x\cdot x}]&=[l_{x,e},l_{x,\bar x}]\\
  &=l_{xex,\bar x}-l_{x,ex\bar x}\\
  &=l_{x\cdot x,\bar x}-l_{x,\bar x\cdot\bar x}\quad\text{by \eqref{eq:G5},}\\
  &=\left( [L_{x\cdot x},L_x]+L_{x\cdot(x\cdot x)}\right) -
    \left([L_x,L_{x\cdot x}]+L_{(x\cdot x)\cdot x}\right)
     \quad\text{by \eqref{eq:G8},}\\
  &=-2[L_x,L_{x\cdot x}]\quad\text{by flexibility.}
\end{split}
\]
Since the characteristic is $\ne 3$, this shows that
$[L_x,L_{x^{\cdot 2}}]=0$ for any $x\in J$ which, together with
the flexible law, shows that $J^{(e)}$ is a noncommutative Jordan
algebra, thus completing the proof of the first part of the
Theorem.

\smallskip

Conversely, let $(A,\cdot)$ be a unital algebra in $\cV$ endowed
with an involution $x\mapsto \bar x$. Let $e=1_A$ be the unity
element and use \eqref{eq:bintriple} to define a triple product on
$A$. Then
\[
\begin{split}
xey&= x\cdot (\bar e\cdot y)-
  \bar e\cdot (x\cdot y)+(\bar e\cdot x)\cdot y\\
   &=x\cdot y-x\cdot y+x\cdot y,
\end{split}
\]
so $x\cdot y=xey$ for any $x,y\in A$, and hence $eee=e$ and
$eex=x=xee$ for any $x\in A$. Also,
\[
exe=e\cdot(\bar x\cdot x)-\bar x\cdot (e\cdot e)
 +(\bar x\cdot e)\cdot e=\bar x
\]
for any $x\in A$. Hence, conditions (i)--(iii) are satisfied.

Finally, with $D_{x,y}$ defined by \eqref{eq:Dxy}, we obtain for
any $x,y,z\in A$,
\begin{equation}\label{eq:G10}
D_{x,y}\bar z=\overline{D_{\bar x,\bar y}z}
 =-\overline{D_{\bar y,\bar x}z},
\end{equation}
because
\[
(x,y,\bar z)^\cdot =-\overline{(z,\bar y,\bar x)^\cdot}=
 \overline{(\bar x,\bar y,z)^\cdot}
\]
by flexibility. Since $D_{x,y}\in \Der A$, \eqref{eq:G10} and
\eqref{eq:bintriple} show that
\begin{equation}\label{eq:G11}
[D_{x,y},l_{a,b}]=l_{D_{x,y}a,b}-l_{a,D_{\bar y,\bar x}b}
\end{equation}
for any $x,y,a,b\in A$. Now, \eqref{eq:cyclic} amounts to
\begin{equation}\label{eq:G12}
\begin{split}
[L_{a\cdot b},L_c]&+[L_{b\cdot c},L_a]+[L_{c\cdot a},L_b]\\
 &= L_{[a\cdot b,c]+[b\cdot c,a]+[c\cdot a,b]}
\end{split}
\end{equation}
for any $a,b,c\in A$. But
\[
\begin{split}
&L_{[a\cdot b,c]+[b\cdot c,a]+[c\cdot a,b]}\\
&\quad =L_{(a,b,c)^\cdot +(b,c,a)^\cdot+(c,a,b)^\cdot}\\
&\quad =L_{(a,b,c)^\cdot-(b,a,c)^\cdot +(b,c,a)^\cdot}\\
&\quad =L_{D_{a,b}c}+L_{(b,c,a)^\cdot}\\
&\quad =[D_{a,b},L_c]+L_{(b,c,a)^\cdot}\\
&\quad =-[[L_a,L_b],L_c]+[L_{[a,b]},L_c]+L_{(b,c,a)^\cdot},
\end{split}
\]
thus \eqref{eq:G12} becomes,
\[
[L_c,[L_a,L_b]+L_{b\cdot a}] =\left([L_{c\cdot
a},L_b]+L_{b\cdot(c\cdot a)}\right)
  -\left( [L_a,L_{b\cdot c}]+L_{(b\cdot c)\cdot a}\right).
\]
Substituting $b$ by $\bar b$ and using \eqref{eq:bintriple}, this
last equation is equivalent to
\begin{equation}\label{eq:G13}
[L_c,l_{a,b}]=l_{c\cdot a,b}-l_{a,\bar c\cdot b}
\end{equation}
for any $a,b,c\in A$. From \eqref{eq:G8}, \eqref{eq:G11} and
\eqref{eq:G13} we obtain
\[
\begin{split}
[l_{x,y},l_{a,b}]&=[-D_{x,\bar y}+L_{x\cdot \bar y},l_{a,b}]\\
 &=-l_{D_{x,\bar y}a,b}+l_{a,D_{y,\bar x}b}+
   l_{(x\cdot\bar y)\cdot a,b}-l_{a,(y\cdot\bar x)\cdot b}\\
 &=l_{l_{x,y}a,b}-l_{a,l_{y,x}b}\\
 &=l_{xya,b}-l_{a,yxb},
\end{split}
\]
thus proving that $A$ is a GJTS with the triple product defined by
\eqref{eq:bintriple}
\end{proof}

\begin{remark}\label{re:GJTS}
Notice that in the proof of the first part of the Theorem above,
the restriction on the characteristic to be $\ne 3$ has only been
used to prove \eqref{eq:G4} and $[L_x,L_{x^{\cdot 2}}]=0$ for any
$x\in J$.
\end{remark}

\medskip

As a particular case, for Jordan triple systems, the following
known result is recovered \cite[1.4, 3.2]{Loos}:

\begin{corollary}\label{JTS}
Let $J$ be a Jordan triple system over a field $F$ of
characteristic $\ne 2,3$ which contains an element $e\in J$ such
that
\begin{equation}\label{eq:exe}
exe=x
\end{equation}
for any $x\in J$. Then, the homotope algebra $J^{(e)}$ (with the
product $x\cdot y=xey$) is a unital Jordan algebra with $1=e$.
Moreover, for any $x,y,z\in J$:
\begin{equation}\label{eq:tripleJTS}
xyz=x\cdot (y\cdot z)-y\cdot (x\cdot z)+(y\cdot x)\cdot z.
\end{equation}
Conversely, if $A$ is a unital Jordan algebra with unity $e$, then
the triple product $xyz$ defined by \eqref{eq:tripleJTS} becomes a
Jordan triple system satisfying \eqref{eq:exe}.
\end{corollary}
\begin{proof}
For a Jordan triple system $J$ satisfying \eqref{eq:exe}, all the
conditions (i)--(iii) in Theorem \ref{GJTS} are automatically
satisfied. Moreover, for any $x\in J$, $\bar x=exe=x$ so that the
involution law $\overline{x\cdot y}=\bar y\cdot \bar x$ is
equivalent to the commutative law $x\cdot y=y\cdot x$.
\end{proof}

\bigskip

\section{$(-1,-1)$-balanced Freudenthal Kantor triple systems}

An algebra $Q$ over the field $F$ is said to be \emph{quadratic}
if it is unital and for any $x\in Q$, $1$, $x$ and $x^2$ are
linearly dependent. Then the set of \emph{vectors} $V=\{ x\in
Q\setminus F1 : x^2\in F1\}$ is a subspace of $Q$ with $Q=F1\oplus
V$ \cite{Osborn}. For any $u,v\in V$,
\begin{equation}\label{eq:prodvectors}
uv=-(u\vert v)1+u\times v,
\end{equation}
where $(.\vert .)$ is a bilinear form and $\times$ is an
anticommutative multiplication on $V$. The triple $\Bigl(
V,(.\vert .),\times\Bigr)$ determines the algebra $Q$, so we will
write $Q=Q\Bigl( V,(.\vert .),\times\Bigr)$. Moreover, for any
$x\in Q$,
\begin{equation}\label{eq:quad}
x^2-T(x)x+N(x)1=0
\end{equation}
where $T$ is a linear form and $N$ a quadratic form on $A$, called
the \emph{norm}, given for any $\alpha\in F$ and $u\in V$ by
\begin{equation}\label{eq:TN}
\left\{
\begin{aligned}
T(\alpha 1+u)&=2\alpha\\
N(\alpha 1+u)&=\alpha^2+(u\vert u)
\end{aligned}\right.
\end{equation}
In particular, $T(x)=N(x,1)$ for any $x$, where
$N(x,y)=N(x+y)-N(x)-N(y)$ is the associated symmetric bilinear
form.

It is well-known \cite{Osborn} that the bilinear form $(.\vert .)$
is symmetric if and only if the map $x\mapsto \bar x=T(x)1-x$ is
an involution (the \emph{standard involution}), and that the
quadratic algebra $Q$ is flexible if and only if $(.\vert .)$ is
symmetric and $(u\times v\vert w)=(u\vert v\times w)$ for any
$u,v,w\in V$. Notice that if $(.\vert .)$ is symmetric, it is
determined by $N$ and that any flexible quadratic algebra is a
noncommutative Jordan algebra.

\begin{theorem}\label{grado2}
Let $S$ be a $(-1,-1)$-balanced Freudenthal Kantor triple system
over a field $F$ and let $e\in S$ such that $\langle e\vert
e\rangle \ne 0$ ($\langle .\vert .\rangle$ as in
\eqref{eq:FKTSdos}). Define a binary product on $S$ by
\[
x\cdot y=\frac{1}{\langle e\vert e\rangle}exy
\]
for any $x,y\in S$. Then $(S,\cdot)$ is a quadratic algebra in the
variety $\cV$ with norm given by $N(x)=\frac{\langle x\vert
x\rangle}{\langle e\vert e\rangle}$ for any $x\in S$. Moreover,
the original triple product on $S$ is related to the binary
product by
\begin{equation}\label{eq:binter2}
xyz=\langle e\vert e\rangle\Bigl((\bar x\cdot y)\cdot z-\bar
x\cdot (y\cdot z)+y\cdot(\bar x\cdot z)\Bigr)
\end{equation}
for any $x,y,z\in S$, where $x\mapsto \bar x$ denotes the standard
involution of the quadratic algebra $(S,\cdot)$.

Conversely, let $(Q,\cdot)$ be a quadratic algebra in $\cV$ with
norm $N$ and define a triple product on $Q$ by the formula
\begin{equation}\label{eq:binter3}
xyz=(\bar x\cdot y)\cdot z-\bar x\cdot (y\cdot z)+y\cdot(\bar
x\cdot z),
\end{equation}
where $x\mapsto \bar x=N(x,1)1-x$ is the standard involution on
$Q$. Then $Q$, with this triple product, is a $(-1,-1)$-BFKTS with
associated nonzero symmetric bilinear form given by $\langle
x\vert x\rangle =N(x)$ for any $x\in Q$.
\end{theorem}

\begin{proof}
Let $S$ be a $(-1,-1)$-BFKTS and let $e\in S$ with $\langle e\vert
e\rangle\ne 0$. Define a new triple product on $S$ by the formula
\[
\widetilde{xyz}=\frac{1}{\langle e\vert e\rangle}yxz.
\]
Then the map $\tilde l_{x,y}: z\mapsto \widetilde{xyz}$ equals
\[
\begin{split}
\tilde l_{x,y}&=\frac{1}{\langle e\vert e\rangle}l_{y,x}\\
 &=\frac{1}{\langle e\vert e\rangle}\Bigl(2\langle x\vert
 y\rangle id-l_{x,y}\Bigr),
\end{split}
\]
where $l_{x,y}$ is the `left multiplication' on $S$ (see
\eqref{eq:FKTSdos}). Thus, for any $x,y,z,t\in S$,
\begin{equation}\label{eq:gorro}
\begin{split}
[\tilde l_{x,y},\tilde l_{z,t}]&=
  \frac{1}{\langle e\vert e\rangle^2}[l_{y,x},l_{t,z}]\\
  &=\frac{1}{\langle e\vert e\rangle^2}
     \Bigl(l_{yxt,z}-l_{t,xyz}\Bigr)\\
  &=\frac{1}{\langle e\vert e\rangle^2}
     \Bigl(l_{t,yxz-2\langle x\vert y\rangle z}-
        l_{xyt-2\langle x\vert y\rangle t,z}\Bigr)\\
  &=\frac{1}{\langle e\vert e\rangle^2}
     \Bigl(l_{t,yxz}-l_{xyt,z}\Bigr)\\
  &=\tilde l_{\widetilde{xyz},t}-\tilde l_{z,\widetilde{yxt}}.
\end{split}
\end{equation}
Therefore, $(S,\widetilde{xyz})$ is a GJTS. Moreover,
\[
\left\{
\begin{aligned}
\widetilde{eex}&=\frac{1}{\langle e\vert e\rangle} eex=x,\\
\widetilde{xee}&=\frac{1}{\langle e\vert e\rangle} exe=x,
\end{aligned}\right.
\]
so \eqref{eq:G4} is satisfied, and
\[
\begin{split}
\bar x= \widetilde{exe}=\frac{1}{\langle e\vert e\rangle}xee
  &=\frac{1}{\langle e\vert e\rangle}
        \Bigl(-eex+2\langle e\vert x\rangle e\Bigr)\\
  &=-x+N(e,x)e
\end{split}
\]
with $N$ the quadratic form given by $N(x)=\frac{\langle x\vert
x\rangle }{\langle e\vert e\rangle}$ for any $x\in S$. Define the
algebra $(S,\cdot)$ by means of
\[
x\cdot y=\frac{1}{\langle e\vert e\rangle}exy=\widetilde{xey}
\]
for any $x,y\in S$. This algebra $(S,\cdot)$ is unital with
$1_S=e$ and for any $x\in S$:
\[
\begin{split}
x\cdot x=\frac{1}{\langle e\vert e\rangle}exx
  &=\frac{1}{\langle e\vert e\rangle}
    \Bigl(2\langle e\vert x\rangle-xex\Bigr)\\
  &=N(e,x)x-N(x)e,
\end{split}
\]
so $(S,\cdot)$ is quadratic. Now by Theorem \ref{GJTS} and Remark
\ref{re:GJTS}, the algebra $(S,\cdot)$ is flexible (and hence
noncommutative Jordan) and satisfies that $D_{x,y}\in
\Der(S,\cdot)$ for any $x,y$. Finally, from \eqref{eq:bintriple},
\[
\begin{split}
xyz&=\langle e\vert e\rangle \widetilde{yxz}\\
&=\langle e\vert e\rangle
  \Bigl( y\cdot(\bar x\cdot z)-\bar x\cdot (y\cdot z)
      +(\bar x\cdot y)\cdot z\Bigr)
\end{split}
\]
for any $x,y,z\in S$, thus completing the proof of the first part.

Conversely, if $(Q,\cdot)$ is a quadratic algebra in $\cV$ with
norm $N$ and we use \eqref{eq:binter3} to define a triple product
on $Q$, then
\[
\left\{
\begin{aligned}
xxy&=(\bar x\cdot x)\cdot y-\bar x\cdot(x\cdot y)
   +x\cdot(\bar x\cdot y)=N(x)y\\
xyx&=(\bar x\cdot y)\cdot x-\bar x\cdot (y\cdot x)
   +y\cdot(\bar x \cdot x)=N(x)y
\end{aligned}\right.
\]
since $x\cdot \bar x=\bar x\cdot x=N(x)1$, $\bar x\cdot (x\cdot
y)=x\cdot (\bar x\cdot y)$ and $(\bar x,y,x)^\cdot =0$ for any
$x$, as $\bar x\in F1+Fx$ and $(Q,\cdot)$ is flexible. Also, with
$\widetilde{xyz}=yxz=(\bar y\cdot x)\cdot z-\bar y\cdot (x\cdot
z)+x\cdot (\bar y\cdot z)$, Theorem \ref{GJTS} shows that
$(Q,\widetilde{xyx})$ is a GJTS and, as in \eqref{eq:gorro}, this
shows that so is $(Q,xyz)$, as required.
\end{proof}

We close this section with a result relating the simplicity of a
quadratic algebra in $\cV$ and of the associated $(-1,-1)$-BFKTS,
constructed in the previous Theorem.

\begin{theorem}\label{simplicity}
Let $(Q,\cdot)$ be a quadratic algebra in $\cV$ and let $(Q,xyz)$
be the associated $(-1,-1)$-BFKTS with triple product given by
\eqref{eq:binter3}. Then $(Q,xyz)$ is simple if and only if either
$(Q,\cdot)$ is simple or $(Q,\cdot)$ is isomorphic to the direct
product of two copies of the ground field $F$: $(Q,\cdot)\cong
F\times F$.
\end{theorem}

\begin{proof}
Assume first that $(Q,xyz)$ is simple and let $0\ne I$ be an ideal
of $(Q,\cdot)$. Let $x\mapsto \bar x$ be the standard involution.
Then either $I=\bar I$ and hence, by \eqref{eq:binter3}, $I$ is an
ideal of $(Q,xyz)$, so $I=Q$ by simplicity, or $I\ne \bar I$. In
the latter case, $I+\bar I$ and $I\cap\bar I$ are ideals of
$(Q,\cdot)$ closed under the involution so, by the previous
argument, $I\cap\bar I=0$ and $I+\bar I=Q$. Besides, $I\ne \bar
I$, so there is some element $x\in I$ with $N(x,1)\ne 0$. But for
any $y\in Q$, by \eqref{eq:quad},
\[
x\cdot y+y\cdot x=N(x,1)y+N(y,1)x-N(x,y)1,
\]
so $N(x,1)y-N(x,y)1\in I$ and hence $\{y\in Q :
N(x,y)=0\}\subseteq I$ and the codimension of $I$ (which coincides
with the codimension of $\bar I$) is $1$. The only possibility is
that $\dim I=\dim \bar I=1$ and $Q=I\oplus\bar I$. In particular
both $I$ and $\bar I$ are one-dimensional quotients of the unital
algebra $Q$, so both $I$ and $\bar I$ are isomorphic, as algebras,
to the ground field $F$, and $(Q,\cdot)\cong F\times F$, as
required.

Conversely, any ideal $I$ of $(Q,xyz)$ satisfies that for any
$x\in I$, $\bar x=x11\in I$, so that $I$ is an ideal of
$(Q,\cdot)$ closed under the involution, and hence $I$ is trivial
in both cases: $(Q,\cdot)$ simple or $(Q,\cdot)\cong F\times F$.
\end{proof}

\section{Simple algebras}

According to Theorems \ref{grado2} and \ref{simplicity}, to obtain
the simple finite dimensional quadratic algebras in $\cV$ it is
enough to consider the simple finite dimensional $(-1,-1)$-BFKTS's
$S$ with an element $e\in S$ such that $\langle e\vert e\rangle=1$
and to define the associated quadratic algebras $(S,\cdot)$ where
\begin{equation}\label{eq:binary}
x\cdot y=exy
\end{equation}
for any $x,y\in S$. The element $e$ becomes the unity element of
$(S,\cdot)$.

The classification of the simple finite dimensional
$(-1,-1)$-BFKTS's over fields of characteristic zero was obtained
in \cite{EKO}. Here we will review the list of examples that
appear in \cite[Section 3]{EKO} and obtain the associated
quadratic algebras in $\cV$. This will set the stage for the
classification in the last section.

\smallskip

\noindent\textbf{5.(i) Orthogonal type:}

Let $S$ be a vector space endowed with a symmetric bilinear form
$\langle .\vert .\rangle$ and an element $e\in S$ such that
$\langle e\vert e\rangle =1$. Then $S$ becomes a $(-1,-1)$-BFKTS
with the triple product
\[
xyz=\langle z\vert x\rangle y
  -\langle z\vert y\rangle x+\langle x\vert y\rangle z
\]
for any $x,y,z\in S$. Therefore, \eqref{eq:binary} becomes
\[
x\cdot y=\langle e\vert y\rangle x +\langle e\vert x\rangle y
  -\langle x\vert y\rangle e
\]
so that $(S,\cdot)=Fe\oplus V$, with $V=(Fe)^\perp$ (the
orthogonal subspace to $Fe$ relative to $\langle .\vert
.\rangle$), is the Jordan algebra of a quadratic form: for any
$u,v\in V$ and $\alpha,\beta\in F$,
\begin{equation}\label{eq:orth}
(\alpha e+u)\cdot (\beta e+v)=
 (\alpha\beta-\langle u\vert v\rangle)e+(\alpha v+\beta u).
\end{equation}

\smallskip

\noindent\textbf{5.(ii) Unitarian type:}

Let $K$ be a quadratic \'etale $F$-algebra; that is, either $K$ is
a quadratic field extension of $F$ (recall that the characteristic
of $F$ is not $2$) or it is isomorphic to $F\times F$; and let $S$
be a left $K$-module endowed with a hermitian form $h: S\times
S\rightarrow K$ and an element $e\in S$ with $h(e,e)=1$. Thus, $h$
is $F$-bilinear and
\[
\begin{split}
h(\alpha x,y)&=\alpha h(x,y)\\
h(x,y)&=\overline{h(y,x)}
\end{split}
\]
for any $\alpha\in F$ and $x,y\in S$, where $\alpha\mapsto
\bar\alpha$ is the nontrivial $F$-automorphism of $K$.

Then $S$ is a $(-1,-1)$-BFKTS with the triple product
\[
xyz=h(z,x)y-h(z,y)x+h(x,y)z
\]
for any $x,y,z\in S$. Thus \eqref{eq:binary} becomes here:
\begin{equation}\label{eq:unitarian1}
x\cdot y=h(y,e)x+h(e,x)y-h(y,x)e
\end{equation}
for any $x,y\in S$; so that $(S,\cdot)=Ke\oplus W$, with
$W=(Ke)^\perp=\{x\in S: h(e,x)=0\}$ and for any $\alpha,\beta\in
K$ and $u,v\in W$:
\begin{equation}\label{eq:unitarian}
(\alpha e+u)\cdot (\beta e+v)=
 (\alpha\beta-h(v,u))e+(\bar\alpha v+\beta u).
\end{equation}
Therefore, $(S,\cdot)$ is the \emph{structurable algebra}
associated to the hermitian form $-h\vert_{W}$ (see \cite[\S\, 8,
Example (iii)]{Allison}.

\smallskip

\noindent\textbf{5.(iii) Symplectic type:}

Change $K$ to $H$, a quaternion algebra over $F$, in the unitarian
type; so that now $S$ is a left $H$-module endowed with a
hermitian form $h:S\times S\rightarrow H$ and an element $e\in S$
with $h(e,e)=1$. (Here $\alpha\mapsto \bar \alpha$ denotes the
standard involution in $H$.) As before, $S=He\oplus W$ with
$W=(He)^\perp$, but now \eqref{eq:unitarian} becomes
\begin{equation}\label{symplectic}
(\alpha e+u)\cdot (\beta e+v)=
 \bigl(\bar\alpha\beta+\beta(\alpha-\bar\alpha)-h(v,u)\bigr)e
  +(\bar\alpha v+\beta u).
\end{equation}
for any $\alpha,\beta\in H$ and $u,v\in W$.

\bigskip

In order to deal with the remaining types, some preliminaries are
needed.

Given a quadratic algebra $Q=Q\Bigl(V,(.\vert .),\times\Bigr)$ and
a nonzero scalar $\mu\in F$, we will denote by $Q^{[\mu]}$ the
quadratic algebra
\[
Q^{[\mu]}=Q\Bigl(V,\mu(.\vert .),\times\Bigr).
\]
(Same anticommutative multiplication on $V$, but bilinear forms
scaled by $\mu$.)

There is a related construction in the literature. Given any
algebra $A$ and a scalar $\alpha\in F$, the \emph{scalar mutation}
$A^{(\alpha)}$ (see \cite{Albert,McCrimmon}) is the algebra
defined on the same vector space but with the new product
\[
x\buildrel \alpha\over{\cdot} y=\alpha xy+(1-\alpha)yx
\]
for any $x,y\in A$. For a flexible quadratic algebra
$Q=Q\Bigl(V,(.\vert .),\times\Bigr)$, it follows immediately from
\eqref{eq:prodvectors} that $Q^{(\alpha)}=Q\Bigl(V,(.\vert
.),(2\alpha -1)\times\Bigr)$ (same bilinear form, but
anticommutative multiplication scaled by $2\alpha -1$). Also, for
any $0\ne \nu\in F$, the linear endomorphism of $Q=F1\oplus V$,
given by $\varphi(1)=1$ and $\varphi(v)=\nu v$ for any $v\in V$,
gives an isomorphism $Q^{[\nu^2]}\cong Q\Bigl(V,(.\vert
.),\nu^{-1}\times\Bigr)$, so that for $\alpha\ne \frac{1}{2}$, the
scalar mutation $Q^{(\alpha)}$ is isomorphic to
$Q^{\left[\frac{1}{(2\alpha -1)^2}\right]}$.

\bigskip

\noindent\textbf{5.(iv) $\mathbf{D}_{\boldsymbol{\mu}}$-type:}

Let $S$ be a four dimensional vector space endowed with a
nondegenerate symmetric bilinear form $\langle .\vert .\rangle$
and an element $e$ with $\langle e\vert e\rangle =1$, and let
$\phi: S\times S\times S\times S\rightarrow F$ be a nonzero
alternating multilinear form (unique up to multiplication by a
nonzero scalar). Define the alternating triple product $[xyz]$ on
$S$ by means of
\[
\phi(x,y,z,t)=\langle [xyz]\vert t\rangle
\]
for any $x,y,z,t\in S$. Then \cite[Lemma 3.2]{EKO} there exists a
nonzero scalar $\mu\in F$ such that
\begin{equation}\label{eq:mu}
\langle [a_1a_2a_3]\vert [b_1b_2b_3]\rangle
 =\mu\det\Bigl(\langle a_i\vert b_j\rangle\Bigr)
\end{equation}
for any $a_i,b_i\in S$ ($i=1,2,3$).

In this case, $S$ becomes a $(-1,-1)$-BFKTS with the triple
product
\[
xyz=[xyz]+\langle z\vert x\rangle y-\langle z\vert y\rangle x
  +\langle x\vert y\rangle z
\]
for any $x,y,z\in S$. Thus, \eqref{eq:binary} becomes
\[
x\cdot y=[exy]+\langle e\vert y\rangle x
  +\langle e\vert x\rangle y - \langle x\vert y\rangle e,
\]
so that $(S,\cdot)=Fe\oplus V$ with $V=(Fe)^\perp$ and for any
$\alpha,\beta\in F$ and $u,v\in V$,
\begin{equation}\label{eq:Dmu}
(\alpha e+u)\cdot (\beta e+v) =
 (\alpha\beta -\langle u\vert v\rangle)e
  +(\alpha v+\beta u +u\times v),
\end{equation}
where $u\times v=[euv]$. That is $(S,\cdot)=Q\Bigl(V,\langle
.\vert .\rangle,\times\Bigr)$. From \eqref{eq:mu} and since
$\langle e\vert e\rangle =1$, it follows that for any $u,v\in V$,
\[
\langle u\times v\vert u\times v\rangle =
 \mu\begin{vmatrix}
 \langle u\vert u\rangle &\langle u\vert v\rangle\\
 \langle u\vert v\rangle &\langle v\vert v\rangle
 \end{vmatrix},
\]
so that
\begin{equation}\label{eq:cross}
(u\times v\vert u\times v)=
\begin{vmatrix}
(u\vert u)&(u\vert v)\\ (u\vert v)&(v\vert v)
\end{vmatrix},
\end{equation}
where $(u\vert v)=\mu\langle u\vert v\rangle$ for any $u,v\in V$.
The above equation \eqref{eq:cross} shows that $\times$ is a
vector cross product on $V$ relative to the nondegenerate
symmetric bilinear form $(.\vert .)$ (see \cite{BrownGray}) and,
therefore, the quadratic algebra $H=Q\Bigl(V,(.\vert
.),\times\Bigr)$ is a quaternion algebra over $F$. But then we
conclude that
\[
(S,\cdot)=Q\Bigl(V,\langle.\vert .\rangle,\times\Bigr)
 =Q\Bigl(V,\mu^{-1}(.\vert .),\times\Bigr)=H^{[\mu^{-1}]}.
\]
That is, $(S,\cdot)$ is a quadratic algebra obtained from a
quaternion algebra by scaling the bilinear form on the subspace of
vectors.

Conversely, it is straightforward to check that for any quaternion
algebra $H$ and nonzero scalar $\nu\in F$, the quadratic algebra
$H^{[\nu]}$ belongs to $\cV$.

\smallskip

\noindent\textbf{5.(v) $\mathbf{G}$-type:}

Let $C$ be a Cayley algebra (that is, an eight dimensional unital
composition algebra) over $F$ with norm $n$ and trace $t$ and let
$S=C_0=\{ x\in C : t(x)=0\}$. Let $0\ne \alpha\in F$ and consider
the nondegenerate symmetric bilinear form and the triple product
on $S$ given by:
\[
\left\{
\begin{aligned}
&\langle x\vert y\rangle =-2\alpha t(xy)\\
&xyz=\alpha\Bigl(D_{x,y}(z)-2t(xy)z\Bigr)
\end{aligned}\right.
\]
for any $x,y,z\in S$, where
\[
D_{x,y}=[L_x,L_y]+[L_x,R_y]+[R_x,R_y]
\]
(a derivation of $C$). We refer to \cite[Chapter III]{Schafer} for
the basic properties of Cayley algebras. Assume that there is an
element $e\in S$ with $\langle e\vert e\rangle=1$. Then
$t(e^2)=-2n(e)=-\frac{1}{2\alpha}$, so
\[
\left\{
\begin{aligned}
&\langle x\vert y\rangle = \frac{t(xy)}{t(e^2)},\\
&xyz=\frac{1}{4n(e)}\Bigl(D_{x,y}(z)-2t(xy)z\Bigr) .
\end{aligned}\right.
\]
Here $K=F1+Fe$ is a quadratic \'etale subalgebra of $C$ and
$S=Fe\oplus V$, where $V=\{ x\in C_0 : t(ex)=0\}=\{ x\in C :
t(Kx)=0\}$. For any $u,v\in V$, \eqref{eq:binary} becomes
\[
u\cdot v=\frac{1}{4n(e)}D_{e,u}(v),
\]
but $D_{x,y}(z)=[[x,y],z]+3(x,z,y)$ (\cite[(3.70)]{Schafer}), so
\[
\begin{split}
D_{e,u}(v)&=[[e,u],v]+3(e,v,u)\\
&=-2[ue,v]+3(u,e,v)\quad\text{(since $eu+ue=0$ as $t(eu)=0$)}\\
&=(ue)v+2v(ue)-3u(ev).
\end{split}
\]
Also, $t(eu)=t(ev)=0$ so that, by alternativity, for any $u,v\in
V$,
\[
\begin{split}
(ue)v&=-(uv)e,\\
v(ue)&=-v(eu)=e(vu)=-n(u,v)e-e(uv),\\
u(ev)&=-e(uv).
\end{split}
\]
Hence $D_{e,u}(v)=-2n(u,v)e+[e,uv]$ and \eqref{eq:binary} becomes
\begin{equation}\label{eq:Gtype}
(\alpha e+u)\cdot (\beta e+v) =\Bigl(\alpha\beta
-\frac{n(u,v)}{2n(e)}\Bigr)e +
  \Bigl( \alpha v+\beta u+\frac{1}{4n(e)}[e,uv]\Bigr)
\end{equation}
for any $\alpha,\beta\in F$ and $u,v\in V$.

Now (\cite[\S 6]{Jac58}, \cite[\S\S\, 2,3]{EMy}) for any $u,v\in
V$,
\begin{equation}\label{eq:color}
uv=-\sigma(u,v)+u*v,
\end{equation}
where $\sigma(x,y)=\frac{1}{2}\Bigl(
n(x,y)-\frac{1}{n(e)}n(ex,y)e\Bigr)$ for any $x,y\in V$. Then
$\sigma: V\times V\rightarrow K$ is a hermitian form, $*$ is
anticommutative and
\begin{equation}\label{eq:colorprop}
\left\{
\begin{aligned}
&\mu(x*y)=(\bar \mu x)*y,\\
&\sigma(x,y*z)=\overline{\sigma(x*y,z)},\\
&(x*y)*z=\sigma(x,z)y-\sigma(y,z)x,
\end{aligned}\right.
\end{equation}
for any $\mu\in K$ and $x,y\in V$. The quadratic algebra
$B=F1\oplus V=Q\Bigl(V,-n(.\vert .),*\Bigr)$ is a \emph{color
algebra} (see \cite{EMy}) and any color algebra is obtained in
this way. For the origin and basic properties of color algebras
one may consult \cite{EMy} and the references therein.

\begin{proposition}
The linear map $\varphi: B^{[-2]}\rightarrow (S,\cdot)$ given by
$\varphi(1)=e$ and $\varphi(u)=-2eu$ for any $u\in V$ is an
isomorphism of algebras.
\end{proposition}

\begin{proof}
Since $e$ is the unity element of $(S,\cdot)$ and
$B^{[-2]}=Q\Bigl(V,2n(.\vert .),*\Bigr)$, it is enough to prove
that for any $u,v\in V$,
\begin{equation}\label{eq:phiuphiv}
\varphi\bigl( -2n(u,v)1+u*v)=\varphi(u)\cdot\varphi(v).
\end{equation}
But
\begin{equation}\label{eq:aprobar}
\begin{split}
&\varphi\bigl(-2n(u,v)1+u*v)=-2n(u,v)e-2e(u*v),\\
&\varphi(u)\cdot\varphi(v)=4(eu)\cdot (ev)
 =-2\frac{n(eu,ev)}{n(e)}e+\frac{1}{n(e)}[e,(eu)(ev)]
\end{split}
\end{equation}
and $n(eu,ev)=n(e)n(u,v)$ by the composition property of the norm
of $C$, while \eqref{eq:color} and \eqref{eq:colorprop} give
\[
\begin{split}
(eu)(ev)&=-\sigma(eu,ev)+(eu)*(ev)\\
 &=-e\bar e\sigma(u,v)+{\bar e}^2(u*v)\\
 &=-n(e)\bigl(\sigma(u,v)+u*v\bigr).
\end{split}
\]
(Note that $\bar e=-e$ and $e\bar e=n(e)1$.) Since $[e,K]=0$ and
$eu+ue=-n(e,u)1=0$, for any $u,v\in V$:
\[
[e,(eu)(ev)]=-n(e)[e,u*v]=-2n(e)e(u*v).
\]
This, together with \eqref{eq:aprobar} proves the validity of
\eqref{eq:phiuphiv}.
\end{proof}

Therefore, the quadratic algebras in $\cV$ associated to the
$(-1,-1)$-BFKTS's of $G$-type are precisely the algebras
$B^{[-2]}$, where $B$ is a color algebra.

\smallskip

\noindent\textbf{5.(vi) $F$-type:}

Here the characteristic of $F$ will be assumed to be $\ne 2,3$.
Let $S$ be an eight dimensional vector space over $F$ endowed with
a nondegenerate symmetric bilinear form $\langle .\vert .\rangle$,
an element $e\in s$ with $\langle e\vert e\rangle =1$ and a
$3$-fold vector cross product $X$ of type I associated to $\langle
.\vert .\rangle$ (see \cite{Eldternary}, \cite[Ch.~8]{Okubobook}
and the references therein). Then $S$ is a $(-1,-1)$-BFKTS with
the triple product
\[
xyz=\frac{1}{3} X(x,y,z)+\langle z\vert x\rangle y
  -\langle z\vert y\rangle x+\langle x\vert y\rangle z
\]
for any $x,y,z\in S$.

Then $S$ has the structure of a Cayley algebra, denoted by $C$,
with unity element $e$, norm $n(x)=\langle x\vert x\rangle$ and
standard involution $x\mapsto \bar x$, such that
\[
X(x,y,z)=(x\bar y)z+\langle x\vert z\rangle y
  -\langle y\vert z\rangle x -\langle x\vert y\rangle z,
\]
so the triple product above becomes
\[
xyz=\frac{1}{3}\left((x\bar y)z+4\langle x\vert z\rangle y
  -4\langle y\vert z\rangle x +2\langle x\vert y\rangle z\right),
\]
for any $x,y,z\in S$, while \eqref{eq:binary} becomes
\[
\begin{split}
x\cdot y&=
 \frac{1}{3}\Bigl(\bar x y +4\langle e\vert y\rangle x +
     2\langle e\vert x\rangle y
        -4\langle x\vert y\rangle e\Bigr)\\
 &=\frac{1}{3}\Bigl(-xy+4\langle e\vert y\rangle x
   +4\langle e\vert x\rangle y- 4\langle x\vert y\rangle e\Bigr)
\end{split}
\]
for any $x,y\in S$, since $x+\bar x=2\langle e\vert x\rangle e$.

Now, $(S,\cdot)=Fe\oplus V$ with $V=(Fe)^\perp$ (orthogonal
relative to $\langle .\vert .\rangle$ and for any $u,v\in V$,
$uv=-\langle u\vert v\rangle e+\frac{1}{2}[u,v]$. Hence, for any
$\alpha,\beta\in F$ and $u,v\in V$
\[
(\alpha e+u)\cdot (\beta e+v)=
 (\alpha\beta -\langle u\vert v\rangle)e
   -\frac{1}{3}\left(\frac{1}{2}[u,v]\right),
\]
so that $(S,\cdot)=C^{\left(-\frac{1}{3}\right)}$ (scalar
mutation), which is isomorphic to $C^{[9]}$.

Therefore, the quadratic algebras in $\cV$ associated to the
$(-1,-1)$-BFKTS's of $F$-type are precisely the algebras
$C^{[9]}$, where $C$ is a Cayley algebra.

\bigskip

\section{Simple quadratic $\cV$-algebras}

The classification of the simple finite dimensional
$(-1,-1)$-BFKTS's in \cite{EKO} over fields of characteristic $0$,
together with the previous sections, does almost all the work
needed to prove our last result:

\bigskip

\begin{theorem}\label{classification}
Let $(Q,\cdot)$ be a finite dimensional simple quadratic algebra
in the variety $\cV$ over a field $F$ of characteristic $0$. Then,
up to isomorphism, either:
\begin{romanenumerate}
\item
$(Q,\cdot)$ is the Jordan algebra of a nondegenerate quadratic
form, with the exception of $(Q,\cdot)\cong F\times F$.
\item
There exists a quadratic \'etale algebra $K$ over $F$ such that
$Q$ is a free $K$-module of rank $\geq 3$, endowed with a
nondegenerate hermitian form $h:Q\times Q\rightarrow K$ such that
$h(1,1)=1$ and for any $x,y\in Q$,
\begin{equation}\label{eq:iiyiii}
x\cdot y=h(y,1)x-h(y,x)1+h(1,x)y.
\end{equation}
In this case, $(Q,\cdot)$ is the structurable algebra of the
restriction of the nondegenerate hermitian form $-h$ to $\{x\in Q
: h(1,x)=0\}$.
\item
There exists a quaternion algebra $H$ over $F$ such that $Q$ is a
free left $Q$-module of rank $\geq 2$, endowed with a hermitian
form $h:Q\times Q\rightarrow H$ such that \eqref{eq:iiyiii} holds.
\item
There exists a quaternion algebra $H$ over $F$ and a nonzero
scalar $\mu\in F$ such that $(Q,\cdot)=H^{[\mu]}$.
\item
There exists a color algebra $B$ over $F$ such that
$(Q,\cdot)=B^{[-2]}$.
\item
There exists a Cayley algebra $C$ over $F$ such that
$(Q,\cdot)=C^{[9]}$.
\end{romanenumerate}
\bigskip

Moreover, two algebras in different items above are not isomorphic
and:
\begin{romanprimeenumerate}
\item
Two algebras of type \textup{(i)} are isomorphic if and only if
their quadratic forms are isometric.
\item
Two algebras $Q_1$ and $Q_2$ in item \textup{(ii)} with associated
\'etale algebras $K_1$ and $K_2$ and hermitian forms $h_1$ and
$h_2$ are isomorphic if and only if the hermitian pairs
$(Q_1,h_1)$ and $(Q_2,h_2)$ are isomorphic. (That is, there is an
isomorphism of $F$-algebras $\sigma :K_1\rightarrow K_2$ and an
$F$-linear bijection $\varphi: Q_1\rightarrow Q_2$ such that
$h_2\bigl(\varphi(x),\varphi(y)\bigr)=\sigma\bigl(h_1(x,y)\bigr)$
for any $x,y\in Q_1$.)
\item
Two algebras $Q_1$ and $Q_2$ in item \textup{(iii)} with
associated quaternion algebras $H_1$ and $H_2$ and hermitian forms
$h_1$ and $h_2$ are isomorphic if and only if the hermitian pairs
$(Q_1,h_1)$ and $(Q_2,h_2)$ are isomorphic.
\item
Two algebras $H_1^{[\mu_1]}$ and $H_2^{[\mu_2]}$ in item
\textup{(iv)} are isomorphic if and only if so are the quaternion
algebras $H_1$ and $H_2$ and $\mu_1=\mu_2$.
\item
Two algebras $B_1^{[-2]}$ and $B_2^{[-2]}$ in item \textup{(v)}
are isomorphic if and only if so are the color algebras $B_1$ and
$B_2$.
\item
Two algebras $C_1^{[9]}$ and $C_2^{[9]}$ in item \textup{(vi)} are
isomorphic if and only if so are the Cayley algebras $C_1$ and
$C_2$.
\end{romanprimeenumerate}
\end{theorem}

\begin{proof}
That the finite dimensional simple quadratic algebras in $\cV$
over $F$ are precisely the algebras in the assertion of the
Theorem follows directly from Theorem \ref{grado2}, Theorem
\ref{simplicity} and the classification of the simple
$(-1,-1)$-BFKTS's in \cite[Theorem 4.3]{EKO}.

For the isomorphism problem, notice that any isomorphism $\varphi$
between two flexible quadratic algebras satisfies $\varphi(\bar
x)=\overline{\varphi(x)}$ for any $x$, where $x\mapsto \bar x$
denotes the standard involution) and hence extends to an
isomorphism between the corresponding $(-1,-1)$-BFKTS's, because
of \eqref{eq:binter3}. Conversely, any isomorphism between the
corresponding $(-1,-1)$-BFKTS's that matches the unity elements of
the quadratic algebras is indeed an isomorphism of the quadratic
algebras.

Now, the assertions in (i')--(iii') would follow from the
corresponding assertions in \cite[Theorem 4.3]{EKO} if it could be
proved that if the hermitian or quadratic pairs $(Q_1,h_1)$ and
$(Q_2,h_2)$ are isomorphic, the isomorphism can be taken to match
the unity elements. This is a direct consequence of Witt's Theorem
\cite{Scharlau} in case (i') and in cases (ii') and (iii') with
$K$ and $H$ being division algebras. For the split cases in (ii')
($K=F\times F$) and (iii') ($H=\Mat_2(F)$) an extra argument is
needed. First, if $K=F\times F$ and $Q$ is a free $K$-module of
rank $\geq 3$ endowed with a nondegenerate hermitian form $h$,
then (see \cite[p.~358]{EKO}, up to isomorphism, $Q=W\times W^*$
for a vector space $W$ ($W^*$ being its dual) and
$h\bigl((u,f),(v,g)\bigr)=\bigl(g(u),f(v)\bigr)\in F\times F=K$.
Now, given any two elements $(u,f)$ and $(u',f')$ with
$h\bigl((u,f),(u,f)\bigr)=1=h\bigl((u',f'),(u',f')\bigr)$ (that
is, $f(u)=1=f'(u')$), there is a linear bijection
$\varphi:W\rightarrow W$ such that $\varphi(u)=u'$ and
$f\circ\varphi^{-1}=f'$ (just complete $\{u,f\}$ and $\{u',f'\}$
to a couple of dual bases of $W$ and $W^*$). Then the linear map
$\psi: Q\rightarrow Q$ given by
$\psi\bigl((v,g)\bigr)=\bigl(\varphi(v),g\circ\varphi^{-1})$ is an
automorphism of the hermitian pair $(Q,h)$ that carries $(u,f)$ to
$(u',f')$. This finishes the proof of (ii'). Also, if
$H=\Mat_2(F)$ and $Q$ is a free left $H$-module of rank $\geq 2$
endowed with a nondegenerate hermitian form $h$, then
(\cite[p.~359]{EKO}), up to isomorphism, $Q=U\otimes_FW$, where
$U$ is the irreducible (two dimensional) left $H$-module and $W$
is a vector space, and $h(u_1\otimes w_1,u_2\otimes
w_2)=\psi(w_1,w_2)\varphi(-,u_2)u_1\in\End_F(U)=H$, where
$\varphi:U\times U\rightarrow F$ and $\psi:W\times W\rightarrow F$
are nondegenerate skew symmetric bilinear forms. Let $\{a_1,a_2\}$
be a basis of $U$ with $\varphi(a_1,a_2)=1$. Then for any element
$x=a_1\otimes w_1+a_2\otimes w_2\in Q$,
\[
\begin{split}
h(x,x)
 &=\psi(w_1,w_2)\varphi(-,a_2)a_1+\psi(w_2,w_1)\varphi(-,a_1)a_2\\
 &=\psi(w_1,w_2)\Bigl(\varphi(-,a_2)a_1-\varphi(-,a_1)a_2\Bigr)\\
 &=\psi(w_1,w_2)1.
\end{split}
\]
Thus $h(x,x)=1$ if and only if $\psi(w_1,w_2)=1$. But if
$x=a_1\otimes w_1+a_2\otimes w_2$ and $y=a_1\otimes
w_1'+a_2\otimes w_2'$ are two elements of $Q$ with
$h(x,x)=1=h(y,y)$, that is, $\psi(w_1,w_2)=1=\psi(w_1',w_2')$,
there is an element in the symplectic group $\Phi\in Sp(W,\psi)$
such that $\Phi(w_1)=w_1'$ and $\Phi(w_2)=w_2'$. Then
$1\otimes\Phi$ is an automorphism of the hermitian pair $(Q,h)$
that carries $x$ to $y$. This completes the proof of (iii').

For (iv')--(vi') notice that any isomorphism of quadratic algebras
$\varphi:Q\Bigl(V_1,(.\vert .)_1,\times_1\Bigr)\rightarrow
Q\Bigl(V_2,(.\vert .)_2,\times_2\Bigr)$ restricts to an
isomorphism of anticommutative algebras
$\varphi\vert_{V_1}:(V_1,\times_1)\rightarrow (V_2,\times_2)$,
which is also an isometry $\bigl(V_1,(.\vert .)_1\bigr)\rightarrow
\bigl(V_2,(.\vert .)_2\bigr)$. But for a quaternion or Cayley
algebra $Q\Bigl(V,(.\vert .),\times\Bigr)$, the bilinear form
$(.\vert .)$ is determined by $\times$ due to the identity:
\begin{equation}\label{eq:quaca}
(x\times y)\times y=(x\vert y)y-(y\vert y)x,
\end{equation}
and the same happens for color algebras due to the identity:
\begin{equation}\label{eq:colo}
((x\times y)\times y)\times y=\frac{1}{2}(y\vert y)x\times y.
\end{equation}
(See \cite{Elduquecolor}.) Thus, for instance, if
\[
\varphi: H_1^{[\mu_1]}=Q\Bigl(V_1,\mu_1(.\vert .)_1,\times_1\Bigr)
 \rightarrow H_2^{[\mu_2]}=Q\Bigl(V_2,\mu_2(.\vert .)_2,\times_2\Bigr)
\]
is an isomorphism of algebras in item (iv), where $V_i$ is the set
of vectors of the quaternion algebra $H_i=Q\Bigl(V_i,(.\vert
.)_i,\times_i\Bigr)$, $i=1,2$, then
$\psi=\varphi\vert_{V_1}:(V_1,\times_1)\rightarrow (V_2,\times_2)$
is an isomorphism of anticommutative algebras and also an isometry
$\psi:\bigl(V_1,\mu_1(.\vert .)_1\bigr)\rightarrow
\bigl(V_2,\mu_2(.\vert .)_2\bigr)$. By \eqref{eq:quaca}, $\psi$ is
an isometry too between $\bigl(V_1,(.\vert .)_1\bigr)$ and
$\bigl(V_2,(.\vert .)_2\bigr)$, so that $\mu_1=\mu_2$ and
$\varphi$ is an isomorphism between the quaternion algebras $H_1$
and $H_2$. The cases (v') and (vi') follow using the same
arguments (but with \eqref{eq:colo} instead of \eqref{eq:quaca}
for case (v')).
\end{proof}

\bigskip

There appears the natural open question of studying the quadratic
simple algebras in $\cV$ over arbitrary fields of characteristic
$\ne 2$ and check if some other kind of algebras appear. This
would lead to a classification of the simple $(-1,-1)$-BFKTS's
over these fields. The known classification in characteristic $0$
depends heavily on the classification of the simple Lie
superalgebras in characteristic $0$, and the restriction on the
characteristic there is essential.



\providecommand{\bysame}{\leavevmode\hbox
to3em{\hrulefill}\thinspace}
\providecommand{\MR}{\relax\ifhmode\unskip\space\fi MR }
\providecommand{\MRhref}[2]{%
  \href{http://www.ams.org/mathscinet-getitem?mr=#1}{#2}
} \providecommand{\href}[2]{#2}

\end{document}